\documentclass[11pt,draft]{article}
\usepackage{CJK}
\usepackage{amsmath,amsfonts,mathrsfs,amssymb}
\usepackage{indentfirst}

\newcommand{\HT}{\CJKfamily{hei}}

\def\y{\begin{eqnarray*}}
\def\ey{\end{eqnarray*}}

\begin{document}

\title{\HT {Weighted bounds for a class of singular integral operators in variable exponent Herz-Morrey spaces}}

\author{{Yanqi Yang$^{1,2}$\thanks{Corresponding author is Yanqi Yang.
 E-mail address:yangyq@nwnu.edu.cn} ,~~Qi Wu$^{1}$}
\\[8pt]
{1.College of Mathematics and Statistics,  Northwest Normal University, }\\
{ Gansu Lanzhou
730070, P. R. China}\\2.Hubei Key Laboratory of Applied Mathematics (Hubei University), \\{Hubei Wuhan
430062, P. R. China}}
\date{}
\maketitle

{\bf Abstract:}{ Let  T be the singular integral operator with variable kernel defined by
$Tf(x)= p.v. \int_{\mathbb{R}^{n}}K(x,x-y)f(y)\mathrm{d}y$
and $D^{\gamma}(0\leq\gamma\leq1)$  be the fractional differentiation operator, where $K(x,z)=\frac{\Omega(x,z')}{|z|^{n}}$, $z'=\frac{z}{|z|},~~z\neq0$. Let $~T^{\ast}~$and $~T^\sharp~$  be  the adjoint of $T$ and the pseudo-adjoint of $T$,  respectively. In this paper, via the expansion of  spherical harmonics and the estimates of the convolution operators $T_{m,j}$, we shall prove some boundedness results for $TD^{\gamma}-D^{\gamma}T$ and $(T^{\ast}-T^{\sharp})D^{\gamma}$ under natural regularity assumptions on the exponent  function on a class of generalized Herz-Morrey spaces with weight and variable exponent, which extend some known results.  Moreover,  various norm characterizations for  the product $T_{1}T_{2}$ and the pseudo-product $T_{1}\circ T_{2}$ are also established.
{\small
\begin{minipage}{\textwidth}
\end{minipage}



}

\section{Introduction and Main Results}

Given a measurable function
$p:~\mathbb{R}^{n}\longrightarrow[1,\infty)$,  $L^{p(\cdot)}(\mathbb{R}^{n})$ denotes the set of measurable functions $f$ on $\mathbb{R}^{n}$ such that for some $\lambda>0$,
$$
\int_{\mathbb{R}^{n}}\left(\frac{|f(x)|}{\lambda}\right)^{p(x)}\mathrm{d}x<\infty.
$$
Equipped with the Luxemburg-Nakano norm
$$
\|f\|_{L^{p(\cdot)}(\mathbb{R}^{n})}=\inf\left\{\lambda>0:\int_{\mathbb{R}^{n}}\left(\frac{|f(x)|}{\lambda}\right)^{p(x)}\mathrm{d}x\leq1\right\},
$$
then $L^{p(\cdot)}(\mathbb{R}^{n})$ becomes a Banach function space.  If  $p(x)=p_{0}$ is constant,  then $L^{p(\cdot)}(\mathbb{R}^{n})$ equals the usual Lebesgue spaces $L^{p_{0}}(\mathbb{R}^{n})$.

The spaces with variable exponent have been widely studied in recent ten years. The  results show that thay are not only the generalized forms of the classical function spaces with invariable exponent, but also there are some new breakthroughs in the research techniques. These new real variable methods help people further understand the function spaces. Lebesgue spaces with variable exponent $L^{p(\cdot)}(\mathbb{R}^{n})$ become one class of important function spaces due to the fundamental paper [1] by Kov\'{o}\u{c}ik and R\'{a}kosn\'{i}k, Due to the seminal paper, it becomes one of the important class function spaces.  The theory  of the function spaces with variable exponent  have been applied in   fluid dynamics, elastlcity dynamics, calculus of variations and differential equations with non-standard growth conditions(for example, see [2-4]). In [5], authors  proved the extrapolation theorem which leads the boundedness of some classical operators including the commutators on $L^{p(\cdot)}(\mathbb{R}^{n})$. Karlovich and Lerner also  obtained the bundedness of the singular integral commutators in [6]. The boundedness of  some typical operators is being studied with keen interest on spaces with variable exponent (see [7-10]).

The classical theory of Muckenhoupt $A_{p}$ weights is a powerful tool in harmonic analysis, for example in the study of boundary value problems for Laplace's equation on Lipschitz domains. Recently, the generalized Muckenhoupt $A_{p(\cdot)}$ weights with variable $p$ have been intensively studied by many authors, for example, see [11, 12]. Cruz-Uribe and Wang extended the theory of Rubio de Francia extrapolation to the weighted variable Lebesgue $L^{p(\cdot)}(\omega)$ in [12].
 In particular, we note that Diening and H\"{a}st\"{o}  proved the equivalence between the $A_{p(\cdot)}$ conditions and obtained the boundedness of the Hardy-Littlewood maximal operator $M$ on weighted $L^{p(\cdot)}(\mathbb{R}^{n})$ spaces in [13].

It is well known that Herz spaces play an important role in harmonic analysis. After the Herz spaces with one exponent $p(\cdot)$ were introduced in [14], the boundedness of some operators and some characterizations of these spaces were studied widely (see [15-17]).  Herz spaces $\dot{K}^{\alpha,p}_{q(\cdot)}(\omega)$ and $K^{\alpha,p}_{q(\cdot)}(\omega)$ with varoable exponent $A_{p(\cdot)}$ and $p(\cdot)$ weights, were introduced first by Izuki and Noi in [18, 19]. Using the basics on Banach function spaces and the Muckenhoupt theory with variable exponent, they established the boundedness of an intrinsic square function and fractional integral operators on such spaces. At the same time, variable exponent Morrey spaces were introduced and studied in [20] in the Euclidean setting and in [21] in the setting of metric measurable spaces. The weighted variable exponent Morrey spaces $\mathcal{M}^{p(\cdot),\varphi}_{\omega}(\mathbb{R}^{n}$) were introduced first  in [22].  Izuki defined a class of variable exponent Herz-Morrey spaces $M\dot{K}^{\alpha,\lambda}_{p,q(\cdot)}(\mathbb{R}^{n})$  and obtained the boundedness of sublinear operators satisfying a proper size condition in  [23], and studied the well known Hardy-Littlewood-Sobolev theorem in [24], which generalized the one obtained by Lu and Yang [25] in the classical Herz spaces.

Motivated by the results mentioned above, in this paper, we will study the boundedness of the singular integrals with variable kernel and fractional differentiation operators on $M\dot{K}^{\alpha,p}_{q(\cdot)}(\omega)$ and also obtain the results in $\dot{K}^{\alpha,p}_{q(\cdot)}(\omega)$, respectively. In order to do this, we need to recall some notations and definitions.

 The singular integral operator with variable kernel is defined by
$$
Tf(x)= p.v. \int_{\mathbb{R}^{n}}K(x,x-y)f(y)\mathrm{d}y,\eqno(1.1)
$$
where $K(x,z)=\frac{\Omega(x,z')}{|z|^{n}}$, $z'=\frac{z}{|z|},~~z\neq0$, and $\Omega(x,z)$ satisfies the following conditions:
$$
\Omega(x,\lambda z)=\Omega(x,z), ~~{\rm for ~~any}~ x, z \in \mathbb{R}^{n} {\rm and}~~\lambda > 0,\eqno(1.2)
$$
$$
\int_{S^{n-1}}|\Omega(x,z')|^{2}d\sigma(z')\leq A ~{\rm ~~and}~ \int_{S^{n-1}}\Omega(x,z^{\prime})\mathrm{d}\sigma(z^{\prime})=0, ~~{\rm for~any}~ x \in \mathbb{R}^{n},\eqno(1.3)
$$
where  $\mathbb{S}^{n-1}$  denote the unit sphere in $\mathbb{R}^{n} (n\geq2)$ with normalized Lebesgue measure $\mathrm{d}\sigma$.
Noticing that as a  kind of the convolution operator of harmonic analysis, the research about $T$ defined in (1.1) has been widely concerned(see [26-28]),  and  played an important role  in the theory of non-divergent elliptic equations with  discontinuous coefficients(see [29-30]). In the Mihlin conditions, Calder\'{o}n and Zygmund proved the boundness of $T$ on the $L^{2}(\mathbb{R}^{n})$ in [31].

Denote by $D$ the square root of the Laplacian operator satisfing the condition $\widehat{Df}(\xi)=|\xi|\widehat{f}(\xi)$.  Let $T_1$ and $T_2$ be the operators defined in (1.1) which are differentiateded by their kernels $\Omega_1(x,z)$ and $\Omega_2(x,z)$. Let $T_{1}T_{2}$,~$T_{1}\circ T_{2}$ denote the product and pseudo-product of $T_{1}$ and $T_{2}$, respectively.  In 1957,  Calder\'{o}n and Zygmund found that these operators are closely related to the second order linear elliptic equations with variable coefficients and established the following boundedness of the operators $T_{1}^{\ast}$, $T_{1}^{\sharp}$, $T_{1}T_{2}$, $T_{1}\circ T_{2}$ and $ ~D $ on $L^{p}(\mathbb{R}^{n})(1<p<\infty)$(see [32]).\\

\textbf{Theorem A}$^{[32]}$ Let $ 1<p<\infty$,$~\Omega_{1}(x,y) ,~\Omega_{2}(x,y)\in C^{\beta}(C^{\infty})$,$~\beta >1$ satisfy (1.2) and (1.3). Then there is a constant $C$ such that

(1) $\|(T_{1}D-DT_{1})f\|_{L^{p}}\lesssim\|f\|_{L^{p}};$

(2) $\|(T_{1}^{\ast}-T_{1}^{\sharp})Df\|_{L^{p}}\lesssim\|f\|_{L^{p}};$

(3) $\|(T_{1}\circ T_{2}-T_{1}T_{2})Df\|_{L^{p}}\lesssim\|f\|_{L^{p}}.$

With the futher research, Chen and Zhu proved that Theorem A was also true on Weighted Lebesgue space and morrey space in 2015 (see [33]).  After that, Tao and Yang obtained the  boundedness of those operators on the weighted Morrey-Herz spaces and the weak estimates (see [34-35]).  In the next setting, we want to explore the  boundedness of the singular integrals with variable kernel and fractional differentiations on weighted variable exponent Herz-Morrey spaces, especially, we will prove the outcome of Theorem A  on $M\dot{K}^{\alpha,\lambda}_{q,p(\cdot)}(\omega)$.

Let $0\leq\gamma\leq1$.  For tempered distributions $f\in\mathscr{S}^{\prime}(\mathbb{R}^{n})(n=1,2,\ldots)$,  the fractional differentiation operators $\widehat{D^{\gamma}}$ defined by $\widehat{D^{\gamma}f}(\xi)=|\xi|^{\gamma}\widehat{f}(\xi)$, i.e.,
$$
D^{\gamma}f(x)=(|\xi|^{\gamma}\widehat{f}(\xi))^{\vee}(x).
$$

Let $I_{\gamma}$ be the Riesz potential operator of order $\gamma$  defined on the space of tempered distributions modulo polynomials by setting $\widehat{I_{\gamma}f}(\xi)=|\xi|^{-\gamma}
\widehat{f}(\xi).$  It is easy to see that  a locally integrable function $b\in I_{\gamma}(\mathrm{BMO})(\mathbb{R}^{n})$ if and only if  $D^{\gamma}b\in \mathrm{BMO}(\mathbb{R}^{n})$.  Strichartz (see [36])  showed that, for $\gamma\in(0,1),  I_{\gamma}(\mathrm{BMO})(\mathbb{R}^{n})$ is a space of functions modulo constants which is properly contained in $\mathrm{Lip}_{\gamma}(\mathbb{R}^{n})$.

Denote  $\mathcal{H}_{m}$ to be   the space of spherical harmonical homogeneous polynomials of degree $m$. Let  dim$\mathcal{H}_{m}=d_{m}$ and  $\{Y_{m,j}\}_{j=1}^{d_{m}}$ be an orthonormal system of $\mathcal{H}_{m}$. It is well known that  $\{Y_{m,j}\}_{j=1}^{d_{m}},$
$m=0,1,\ldots,$ is a complete orthonormal system in $L^{2}(S^{n-1})$ ( see [32]).  From (1.3), we can expand the function $\Omega(x,z^{\prime})$ in spherical harmonics
$$
\Omega(x,z^{\prime})=\sum_{m\geq0}\sum^{d_{m}}_{j=1}a_{m,j}(x)Y_{m,j}(z^{\prime}),
$$
where
$$
a_{m,j}(x)=\int_{S^{n-1}}\Omega(x,z^{\prime})\overline{Y_{m,j}(z^{\prime})}\mathrm{d}\sigma(z^{\prime}).
$$
Noting that  $\int_{S^{n-1}}\Omega(x,z^{\prime})\mathrm{d}\sigma(z^{\prime})=0$, then $~a_{0,j}= 0$ for any $x \in \mathbb{R}^{n}$.
Let ~$T_{m,j}f(x)=\frac{Y_{m,j}}{|\cdot|^{n}}\ast f(x)$. Then  $T$, defined in (1.1), can be written as
$$
Tf(x)=\sum_{m\geq 1}\sum_{j=1}^{d_{m}}a_{m,j}(x)T_{m,j}f(x).
$$
Let ~$T^{\ast}$ and ~$T^{\sharp}$  be the adjoint of $T$ and the pseudo-adjoint of $T$ respectively, defined by
$$
T^{\ast}f(x)=\sum^{\infty}_{m=1}\sum^{d_{m}}_{j=1}(-1)^{m}T_{m,j}(\overline{a}_{m,j}f)(x)
$$
and
$$
T^{\sharp}f(x)=\sum^{\infty}_{m=1}\sum^{d_{m}}_{j=1}(-1)^{m}\overline{a}_{m,j}(x)T_{m,j}f(x).
$$
We now formulate our other main results as follows.\\
$\textbf{{Theorem}~1.1}$ Let $0\leq\lambda<\infty,~p(\cdot)\in \mathcal{P}(\mathbb{R}^{n})\cap  LH(\mathbb{R}^{n}),~0<q<\infty,~\omega\in A_{rp(\cdot)}$, $1/p_{-}<r<1$  and $-n\delta+\lambda<\alpha<n(1-r)+\lambda$, where $0<\delta<1$ is the constant appearing in Lemma 2.2. Let $\Omega(x,y)$ which satisfies(1.2), (1.3) meet a condition
$$
\max_{|j|\leq2n}\|D_{x}^{\gamma}(\partial^{j}/{\partial y}^{j})\Omega(x,y)\|_{L^{\infty}(\mathbb{R}^{n}\times S^{n-1})}<\infty. \eqno(1.4)
$$
Then we have

(1) $\|(TD^{\gamma}-D^{\gamma}T)f\|_{M\dot{K}^{\alpha,\lambda}_{q,p(\cdot)}(\omega)}\lesssim  \|f\|_{M\dot{K}^{\alpha,\lambda}_{q,p(\cdot)}(\omega)};$

(2) $\|(T^{\ast}-T^{\sharp})D^{\gamma}f\|_{M\dot{K}^{\alpha,\lambda}_{q,p(\cdot)}(\omega)}\lesssim \|f\|_{M\dot{K}^{\alpha,\lambda}_{q,p(\cdot)}(\omega)}.$\\
$\textbf{{Theorem}~1.2}$ Let $0\leq\lambda<\infty,~p(\cdot)\in \mathcal{P}(\mathbb{R}^{n})\cap  LH(\mathbb{R}^{n}),~0<q<\infty,~\omega\in A_{rp(\cdot)}$, $1/p_{-}<r<1$  and $-n\delta+\lambda<\alpha<n(1-r)+\lambda$, where $0<\delta<1$ is the constant appearing in Lemma 2.2. Suppose that $\Omega_{1}(x,y)$ and $\Omega_{2}(x,y)$  satisfy (1.2) and (1.3). If $\Omega_{2}(x,y)$  satisfies (1.4) and
$$
\max_{|j|\leq2n}\|(\partial^{j}/{\partial y}^{j})\Omega_{1}(x,y)\|_{L^{\infty}(\mathbb{R}^{n}\times S^{n-1})}<\infty. \eqno(1.5)
$$
Then we have
$$
\|(T_{1}\circ T_{2}-T_{1}T_{2})D^{\gamma}f\|_{M\dot{K}^{\alpha,\lambda}_{q,p(\cdot)}(\omega)}\lesssim\|f\|_{M\dot{K}^{\alpha,\lambda}_{q,p(\cdot)}(\omega)}.
$$

Next, in the endpoint   case of $\gamma=1$  and $\gamma=0$, $D$ is the square root of Laplacian operator and $D^{0}$ is the identity operator $\mathcal{I}$.  we have the following results.
\\
$\textbf{{Theorem}~1.3}$ Let $0\leq\lambda<\infty,~p(\cdot)\in \mathcal{P}(\mathbb{R}^{n})\cap  LH(\mathbb{R}^{n}),~0<q<\infty,~\omega\in A_{rp(\cdot)}$, $1/p_{-}<r<1$  and $-n\delta+\lambda<\alpha<n(1-r)+\lambda$, where $0<\delta<1$ is the constant appearing in Lemma 2.2.  suppose that $\Omega_i(x,y)(i=1,2)$ satisfies (1.2),(1.3) and (1.5). Then we have

(1) $\|(T_{1}\mathcal{I}-\mathcal{I}T_{1})f\|_{M\dot{K}^{\alpha,\lambda}_{q,p(\cdot)}(\omega)}\lesssim  \|f\|_{M\dot{K}^{\alpha,\lambda}_{q,p(\cdot)}(\omega)};$

(2) $\|(T_{1}^{\ast}-T_{2}^{\sharp})\mathcal{I}f\|_{M\dot{K}^{\alpha,\lambda}_{q,p(\cdot)}(\omega)}\lesssim \|f\|_{M\dot{K}^{\alpha,\lambda}_{q,p(\cdot)}(\omega)};$

(3) $\|(T_{1}\circ T_{2}-T_{1}T_{2})\mathcal{I}f\|_{M\dot{K}^{\alpha,\lambda}_{q,p(\cdot)}(\omega)}\lesssim\|f\|_{M\dot{K}^{\alpha,\lambda}_{q,p(\cdot)}(\omega)}.$\\
$\textbf{{Theorem}~1.4}$ Let $0\leq\lambda<\infty,~p(\cdot)\in \mathcal{P}(\mathbb{R}^{n})\cap  LH(\mathbb{R}^{n}),~0<q<\infty,~\omega\in A_{rp(\cdot)}$, $1/p_{-}<r<1$  and $-n\delta+\lambda<\alpha<n(1-r)+\lambda$, where $0<\delta<1$ is the constant appearing in Lemma 2.2.  suppose that $\Omega(x,y)$ satisfies (1.2), (1.3) and
$$
\max_{|j|\leq2n}\|\nabla_{x}(\partial^{j}/{\partial y}^{j})\Omega(x,y)\|_{L^{\infty}(\mathbb{R}^{n}\times S^{n-1})}<\infty. \eqno(1.6)
$$
Then we have

(1) $\|(TD-DT)f\|_{M\dot{K}^{\alpha,\lambda}_{q,p(\cdot)}(\omega)}\lesssim  \|f\|_{M\dot{K}^{\alpha,\lambda}_{q,p(\cdot)}(\omega)};$

(2) $\|(T^{\ast}-T^{\sharp})Df\|_{M\dot{K}^{\alpha,\lambda}_{q,p(\cdot)}(\omega)}\lesssim \|f\|_{M\dot{K}^{\alpha,\lambda}_{q,p(\cdot)}(\omega)}.$\\
$\textbf{{Theorem}~1.5}$ Let $0\leq\lambda<\infty,~p(\cdot)\in \mathcal{P}(\mathbb{R}^{n})\cap  LH(\mathbb{R}^{n}),~0<q<\infty,~\omega\in A_{rp(\cdot)}$, $1/p_{-}<r<1$  and $-n\delta+\lambda<\alpha<n(1-r)+\lambda$, where $0<\delta<1$ is the constant appearing in Lemma 2.2.  Suppose that $\Omega_{1}(x,y)$ and $\Omega_{2}(x,y)$ satisfies(1.2), (1.3). If $\Omega_{1}(x,y)$ satisfies (1.5) and $\Omega_{2}(x,y)$ satisfies (1.6). Then we have
$$
\|(T_{1}\circ T_{2}-T_{1}T_{2})Df\|_{M\dot{K}^{\alpha,\lambda}_{q,p(\cdot)}(\omega)}\lesssim\|f\|_{M\dot{K}^{\alpha,\lambda}_{q,p(\cdot)}(\omega)}.
$$

We make some conventions. In what follows, $C$ always denotes  a positive constant which is independent of the main parameters involved but whose value may differ in different occurrence. We use the symbol $A\lesssim B$ to denote that there exists a positive constant $C$ such that $A\leq CB$. Denote $P(\mathbb{R}^{n})$ to be the set of $p(\cdot):\mathbb{R}^{n}\rightarrow [1,\infty)$ such that
$$
p^{-}=\mathrm{ess}\inf\{p(x):x\in \mathbb{R}^{n}\}>1,~~p^{+}=\mathrm{ess}\sup\{p(x):x\in \mathbb{R}^{n}\}<\infty.
$$
Let $\mathcal{B}(\mathbb{R}^{n})$ be the set of $p(\cdot)\in P(\mathbb{R}^{n})$ such that the Hardy-Littlewood maximal operator $M$  is bounded on $L^{p(\cdot)}(\mathbb{R}^{n})$.
 For $p(x)\in P(\mathbb{R}^{n})$, we use $p^{\prime}(x)$ to denote the dual exponent of $p(x)$, namely,  $\frac{1}{p(x)}+\frac{1}{p^{\prime}(x)}=1.$ $|E|$ denotes the Lebesgue measure of $E\in \mathbb{R}^{n}$ and $\chi _{E}$ its characteristic function.

\section*{\normalsize\HT \ 2~~~Some Notations and Preliminaries}

\quad In this section, the weighted Herz spaces with variable exponent $\dot{K}^{\alpha,q}_{p(\cdot)}(\omega)$ and the weighted Herz-Morrey spaces with variable exponent $M\dot{K}^{\alpha,\lambda}_{q,p(\cdot)}(\omega)$  will be introduced.
And also, some preliminary lemmas and the estimates of the convolution operators $T_{m,j}$ will be given.

$\textbf{{Definition}~2.1}^{[37]}$   Let $p(\cdot)$ be a measurable function with values in $[1, \infty]$,  we say $p(\cdot)\in LH_{0}(\mathbb{R}^{n})$ if there exist $C>0$ such that for any $x,y\in \mathbb{R}^{n}$, $|x-y|<1/2$,
$$
|p(x)-p(y)|\leq \frac{C}{-\log(|x-y|)},\eqno(2.1)
$$
then $p(\cdot)$ is said local $\log$-H\"{o}lder continuous on $\mathbb{R}^{n}$.

We say $p(\cdot)\in LH_{\infty}(\mathbb{R}^{n})$  if there exist $p_{\infty}\in \mathbb{R}$ and a constant $C>0$ such that all $x\in \mathbb{R}^{n}$,
$$
|p(x)-p(\infty)|\leq \frac{C}{\log(e+|x|)},\eqno(2.2)
$$
where $p_{\infty}=\lim\limits_{x\rightarrow\infty}p(x)>1$,
then $p(\cdot)$ is said $\log$-H\"{o}lder continuous at infinity.

The set of $p(\cdot)$ satisfying (2.1) and (2.2) is denoted by $LH(\mathbb{R}^{n})$. We can easily check that  if $p(\cdot)\in \mathcal{P}(\mathbb{R}^{n})\bigcap  LH(\mathbb{R}^{n})$,
then $p^{\prime}(\cdot)\in \mathcal{P}(\mathbb{R}^{n})\bigcap  LH(\mathbb{R}^{n})$. An important consequence of $\log$-H\"{o}lder continuity  is the fact that if $p^{\prime}(\cdot)\in \mathcal{P}(\mathbb{R}^{n})\bigcap  LH(\mathbb{R}^{n})$, then the Hardy-Littlewood maximal operator $M$  defined by
$$
Mf(x)=\sup_{x\in B}\int_{B}|f(y)|\mathrm{d}y,
$$
is bounded on $L^{p(\cdot)}(\mathbb{R}^{n})$ (see [37]).

Let $\omega$ be a weighted function on $\mathbb{R}^{n}$ , that is, $\omega$ is real-valued, non-negative and locally integrable.    $\omega$ is said to be a Muckenhoupt $A_{1}$ weight if
 $$
 M\omega(x)\leq C \omega(x)~~~~a.e. x\in \mathbb{R}^{n}.
 $$
 For $1<p<\infty$, we say that $\omega$ is an $A_{p}$ weight if
 $$
\sup_{B}\left(\frac{1}{|B|}\int_{B}\omega(x)\mathrm{d}x \right) \left(\frac{1}{|B|}\int_{B}\omega(x)^{1-p^{\prime}}\mathrm{d}x\right)^{p-1}<\infty,
 $$

The weighted variable exponent  Lebesgue spaces  is defined by
$$
L^{p(\cdot)}(\omega):=\{f:f\omega^{\frac{1}{p(\cdot)}}\in L^{p(\cdot)}(\mathbb{R}^{n})\}.
$$
Then  $L^{p(\cdot)}(\omega)$ is a Banach space with respect to the norm
 $$
 \|f\|_{L^{p(\cdot)}({\omega})}=\|f\omega^{\frac{1}{p(\cdot)}}\|_{L^{p(\cdot)}(\mathbb{R}^{n})}.
 $$

If $p(\cdot)=p\in (1,\infty)$, then we can easily see that the definition reduces to the classical weighted   Lebesgue spaces $L^{p}(\omega)$ . \\
$\textbf{{Definition}~2.2}^{[18,19]}$  Let $p(\cdot)\in\mathcal{P}(\mathbb{R}^{n})$. A weight $\omega$ is said to be an $A_{p(\cdot)}$ weight if
$$
\sup_{B}\|\omega^{1/p(\cdot)}\chi_{B}\|_{L^{p(\cdot)}(\mathbb{R}^{n})}\|\omega^{-1/p(\cdot)}\chi_{B}\|_{L^{p^{\prime}(\cdot)}(\mathbb{R}^{n})}<\infty.
$$

Cruz-Uribe $et~al$ showed that $\omega\in A_{p(\cdot)}$ if and only if the Hardy-Littlewood maximal operator $M$ is bounded on $L^{p(\cdot)}(\omega)$  [11]. Suppose that $p_{1}(\cdot),~p_{2}(\cdot)\in \mathcal{P}(\mathbb{R}^{n})\bigcap  LH(\mathbb{R}^{n})$  and $p_{1}(\cdot)\leq p_{2}(\cdot)$, it was proved that  $A_{1}\subset A_{p_{1}(\cdot)}\subset A_{p_{2}(\cdot)}$ in [19].

Let $B_{k}=\{x\in \mathbb{R}^{n}:|x|\leq 2^{k}\},~C_{k}=B_{k}\backslash B_{k-1}$ and $\chi _{k}=\chi_{C_{k}}$ be the characteristic function of the set $C_{k}$ for $k\in \mathbb{Z}$.\\
$\textbf{{Definition}~2.3}^{[19]}$  Suppose $p(\cdot)\in\mathcal{P}(\mathbb{R}^{n}), ~0<q\leq\infty$ and $\alpha\in \mathbb{R}$.  The homogeneous weighted Herz space with variable exponent $\dot{K}^{\alpha,q}_{p(\cdot)}(\omega)$ is defined as the set of all $f\in L_{\mathrm{loc}}^{p(\cdot)}(\mathbb{R}^{n}\backslash \{0\}, \omega)$ such that
 $$
\|f\|_{ \dot{K}^{\alpha,q}_{p(\cdot)}(\omega)}:=\left(\sum_{k \in \mathbb{Z}}2^{\alpha kq}\|f\chi_{k}\|^{q}_{L^{p(\cdot)}(\omega)}\right)^{1/q}<\infty,
 $$
with the usual modification when $q=\infty$.\\
 $\textbf{{Definition}~2.4}^{[19]}$  Suppose $p(\cdot)\in\mathcal{P}(\mathbb{R}^{n}), ~0<q\leq\infty$ and $\alpha\in \mathbb{R}$.  The homogeneous weighted Herz spaces with variable exponent $\dot{K}^{\alpha,q}_{p(\cdot)}(\omega)$ is defined as the set of all $f\in L_{\mathrm{loc}}^{p(\cdot)}(\mathbb{R}^{n}\backslash \{0\}, \omega)$ such that
 $$
\|f\|_{ \dot{K}^{\alpha,q}_{p(\cdot)}(\omega)}:=\left(\sum_{k \in \mathbb{Z}}2^{\alpha kq}\|f\chi_{k}\|^{q}_{L^{p(\cdot)}(\omega)}\right)^{1/q}<\infty,
 $$
with the usual modification when $q=\infty$.\\
$\textbf{{Definition}~2.5}^{[19]}$  Suppose $0\leq\lambda<\infty, ~p(\cdot)\in\mathcal{P}(\mathbb{R}^{n})$, $0<q\leq\infty$ and $\alpha\in \mathbb{R}$.  The homogeneous weighted Herz-Morrey  spaces with variable exponent $M\dot{K}^{\alpha,\lambda}_{q,p(\cdot)}(\omega)$ is defined as the set of all $f\in L_{\mathrm{loc}}^{p(\cdot)}(\omega)$ such that
 $$
\|f\|_{M\dot{K}^{\alpha,\lambda}_{q,p(\cdot)}(\omega)}:=\sup_{L\in \mathbb{Z}}2^{-L\lambda}\left(\sum^{L}_{k=-\infty}2^{\alpha kq}\|f\chi_{k}\|^{q}_{L^{p(\cdot)}(\omega)}\right)^{1/q}<\infty,
 $$
 with the usual modification as $q=\infty$.\\
$\textbf{{Remark}~2.1}$ It follows that $M\dot{K}^{\alpha,0}_{q,p(\cdot)}(\omega)=\dot{K}^{\alpha,q}_{p(\cdot)}(\omega)$. In the case $p(\cdot)=p$ is a  constant and $\omega=1$, then $\dot{K}^{\alpha,q}_{p(\cdot)}(\omega)$ and $M\dot{K}^{\alpha,\lambda}_{q,p(\cdot)}(\omega)$  coincide  with the classical Herz spaces and Herz-Morrey spaces, see [38] and [39], respectively.

We need the following lemmas.\\
$\textbf{{Lemmas}~2.1}^{[19]}$  (Generalized H\"{o}lder's Inequality) Let $X$ be a Banach space. If  $f\in X$ and $g\in X^{\prime}$, then we have

 $$
 \int_{\mathbb{R}^{n}}|f(x)g(x)|\mathrm{d}x\lesssim \|f\|_{X}\|g\|_{X^{\prime}}.
 $$

Noticing that for $\omega\in A_{p}(\cdot)$, it is easy to see that $(L^{p(\cdot)}(\mathbb{R}^{n},\omega))^{\prime}=(L^{p^{\prime}(\cdot)}(\mathbb{R}^{n},\omega^{-1})$ [19].\\
$\textbf{{Lemmas}~2.2}^{[19]}$ If $p_{\cdot}(\cdot)\in p_{2}(\cdot)\in \mathcal{P}(\mathbb{R}^{n})\bigcap  LH(\mathbb{R}^{n})$ and $\omega\in A_{rp(\cdot)}$, $1/p_{-}<r<1$.  Then there exists a constant $0<\delta<1$ such that for all $k,j\in \mathbb{Z}$,
$$
\frac{\|\chi_{B_{k}}\|_{L^{p}(\omega)}}{\|\chi_{B_{j}}\|_{L^{p}(\omega)}}\lesssim
\begin{cases}

2^{(k-j)n\delta},~~k\leq j;\\

2^{(k-j)nr},~~k>j.

\end{cases}$$
$\textbf{{Lemmas}~2.3}^{[19]}$ If $X$ is a Banach function space on $\mathbb{R}^{n}$, then  for all balls $B\subset \mathbb{R}^{n}$ we have
$$
\frac{1}{|B|}\|\chi_{B}\|_{X}\|\chi_{B}\|_{X^{\prime}}\lesssim 1.
$$
$\textbf{{Lemmas}~2.4}^{[33]}$ Let $1<p<\infty$ and $\omega\in A_{p}$. Then the $T_{m,j}$ is bounded on $L^{p(\cdot)}(\omega)$, that is
$$
\|T_{m,j}f\|_{L^{p}(\omega)}\lesssim m^{n/2} \|f\|_{L^{p}(\omega)}.
$$

The following important result  had recently been proved by authors in [12].\\
$\textbf{{Lemmas}~2.5}^{[12]}$ Suppose that there exists a constant $1<p_{0}<\infty$ such that for every $\omega_{0}\in A_{p_{0}}$, the inequality
$$
\|f\|_{L^{p_{0}}(\omega_{0})}\lesssim\|g\|_{L^{p_{0}}(\omega_{0})}
$$
holds for all $f \in L^{p_{0}}(\omega_{0})$ and all measurable function $g$. Let $p(\cdot)\in \mathcal{P}(\mathbb{R}^{n})$ and $\omega$ be a weight. If the Hardy-Littlewood maximal operator $M$ is bounded on ${L^{p^{(\cdot)}}(\omega)}$ and on $L^{p^{\prime}(\cdot)}(\omega^{-\frac{1}{p(\cdot)-1}})$, then we have the inequality
$$
\|f\|_{L^{p(\cdot)}(\omega)}\lesssim\|g\|_{L^{p(\cdot)}(\omega)}
$$
holds for all $f\in L^{p(\cdot)}(\omega)$ and all measurable functions $g$.\\

Combing the Lemma 2.4 and 2.5 above, we obtain the following boundedness of the $T_{m,j}$ on weighted variable exponent Lebesgue spaces.\\
$\textbf{{Lemmas}~2.6}$ Let $p(\cdot)\in \mathcal{P}(\mathbb{R}^{n})\bigcap  LH(\mathbb{R}^{n})$ and $\omega\in A_{p}(\cdot)$. Then the $T_{m,j}$ is bounded on $L^{p(\cdot)}(\omega)$, that is
$$
\|T_{m,j}f\|_{L^{p(\cdot)}(\omega)}\lesssim m^{n/2} \|f\|_{L^{p(\cdot)}(\omega)}.
$$\\
$\textbf{{Lemmas}~2.7}^{[19]}$  Suppose $p(\cdot)\in \mathcal{P}(\mathbb{R}^{n})\bigcap  LH(\mathbb{R}^{n})$, then the assumption $\omega\in A_{p}(\cdot)$ implies the equivalence of the following two conditions

(a) M is bounded on $L^{p(\cdot)}(\omega)$;

(b) M is bounded on $ L^{p^{\prime}(\cdot)}(\omega^{-\frac{1}{p(\cdot)-1}})$.\\
$\textbf{{Lemmas}~2.8}^{[40]}$ Suppose $T$ is a Calder\'{o}n-Zygmund  operator associated to a standard kernel $K$. If $0\leq\lambda<\infty,~p(\cdot)\in \mathcal{P}(\mathbb{R}^{n})\bigcap  LH(\mathbb{R}^{n}),~0<q<\infty,~\omega\in A_{rp(\cdot)}$, $1/p_{-}<r<1$  and $-n\delta+\lambda<\alpha<n(1-r)+\lambda$, where $0<\delta<1$ is the constant appearing in Lemma 2.2, then $T$ is bounded on $M\dot{K}^{\alpha,\lambda}_{q,p(\cdot)}(\omega)$.

$\textbf{{Lemmas}~2.9}$  If $0\leq\lambda<\infty,~p(\cdot)\in \mathcal{P}(\mathbb{R}^{n})\cap  LH(\mathbb{R}^{n}),~0<q<\infty,~\omega\in A_{rp(\cdot)}$, $1/p_{-}<r<1$  and $-n\delta+\lambda<\alpha<n(1-r)+\lambda$, where $0<\delta<1$ is the constant appearing in Lemma 2.2. Then we have
$$
\|T_{m,l}(f)\|_{M\dot{K}^{\alpha,\lambda}_{q,p(\cdot)}(\omega)}\lesssim m^{n/2}\|f\|_{M\dot{K}^{\alpha,\lambda}_{q,p(\cdot)}(\omega)}.
$$

\textbf{Proof of Lemma 2.9}.Let $f \in M\dot{K}^{\alpha,\lambda}_{q,p(\cdot)}(\omega)$, $T_{m,l}f=\frac{Y_{m,l}}{|\cdot|^{n}}\ast f$, so we have the fact that $|Y_{m,l}|\lesssim m^{\frac{n-2}{2}}([32])$. We represent $f$ as
$$
f(x)=\sum_{j=-\infty}^{\infty}f(x)\chi_{j}(x)\ := \sum_{j=-\infty}^{\infty}f_{j}(x).
$$By the definition of the the norm in $M\dot{K}^{\alpha,\lambda}_{q,p(\cdot)}(\omega)$, we have
\begin{align*}
\|T_{m,l}(f)\|^{q}_{M\dot{K}^{\alpha,\lambda}_{q,p(\cdot)}(\omega)}&=\sup_{L\in \mathbb{Z}}2^{-L \lambda q}\sum^{L}_{k=-\infty}2^{\alpha kq}\|T_{m,l}(f)\chi_{k}\|^{q}_{L^{p(\cdot)}(\omega)}\\
&=\sup_{L\in \mathbb{Z}}2^{-L \lambda q}\sum^{L}_{k=-\infty}2^{\alpha kq}\left(\sum_{j=-\infty}^{k-2}\|T_{m,l}(f)\chi_{k}\|^{q}_{L^{p(\cdot)}(\omega)}\right)^{q}\\
&+\sup_{L\in \mathbb{Z}}2^{-L \lambda q}\sum^{L}_{k=-\infty}2^{\alpha kq}\left(\sum_{j=k-1}^{k+1}\|T_{m,l}(f)\chi_{k}\|^{q}_{L^{p(\cdot)}(\omega)}\right)^{q}\\
\end{align*}
\begin{align*}
&+\sup_{L\in \mathbb{Z}}2^{-L \lambda q}\sum^{L}_{k=-\infty}2^{\alpha kq}\left(\sum_{j=k+2}^{\infty}\|T_{m,l}(f)\chi_{k}\|^{q}_{L^{p(\cdot)}(\omega)}\right)^{q}\\
&:=I_{1}+I_{2}+I_{3}.
\end{align*}
For $I_{1}$, we observe that if $x\in C_{k}$, $y\in C_{j}$ and $j\leq k-2$, then $|x-y|\approx |x|\approx 2^{k}$. Then we obtain
$$
 |T_{m,l}f_{j}(x)|
  \lesssim m^{\frac{n-2}{2}}\int_{\mathbb{R}^{n}}\frac{|f_{j}(y)|}{|x-y|^{n}}\mathrm{d}y
 \lesssim m^{\frac{n}{2}}\int_{\mathbb{R}^{n}}\frac{|f_{j}(y)|}{|x-y|^{n}}\mathrm{d}y\approx m^{\frac{n}{2}}\cdot 2^{-kn}\|f_{j}\|_{L^{1}(\mathbb{R}^{n})}.
$$
From above inequality, Lemma 2.1 and Lemma 2.3, we have
\begin{align*}
|T_{m,l}(f_{j})(x)|&\lesssim m^{\frac{n}{2}} \cdot 2^{-kn}\|f\|_{L^{1}(\mathbb{R}^{n})}\\
&\lesssim m^{\frac{n}{2}}\cdot2^{-kn}\|f\omega^{\frac{1}{p(\cdot)}}\chi_{j}\|_{L^{p(\cdot)}(\mathbb{R}^{n})}\|f\omega^{-\frac{1}{p(\cdot)}}\chi_{j}\|_{L^{p^{\prime}(\cdot)}(\mathbb{R}^{n})}\\
&\lesssim m^{\frac{n}{2}} \cdot2^{-kn}\|f\omega^{\frac{1}{p(\cdot)}}\chi_{j}\|_{L^{p(\cdot)}(\mathbb{R}^{n})}\|f\omega^{-\frac{1}{p(\cdot)}}\chi_{B_{j}}\|_{L^{p^{\prime}(\cdot)}(\mathbb{R}^{n})}\\
&\lesssim m^{\frac{n}{2}} \cdot 2^{-kn}\|f\chi_{j}\|_{L^{p(\cdot)}(\omega)}\|\chi_{B_{j}}\|_{(L^{p(\cdot)}(\omega))^{\prime}}\\
&\lesssim m^{\frac{n}{2}} \cdot 2^{(j-k)n}\|f\chi_{j}\|_{L^{p(\cdot)}(\omega)}\|\chi_{B_{j}}\|_{L^{p(\cdot)}(\omega)}^{-1},
\end{align*}
which combing with Lemma 2.2 yields
\begin{align*}
\|T_{m,l}(f_{j})\chi_{k}\|_{L^{p(\cdot)}(\omega)}&\lesssim m^{\frac{n}{2}} \cdot2^{(j-k)n}\|f\chi_{j}\|_{L^{p(\cdot)}(\omega)}\|\chi_{B_{j}}\|_{L^{p(\cdot)}(\omega)}^{-1}\|\chi_{B_{k}}\|_{L^{p(\cdot)}(\omega)}\\
&\lesssim m^{\frac{n}{2}} \cdot 2^{(j-k)n(1-r)}\|f\chi_{j}\|_{L^{p(\cdot)}(\omega)}.\quad\quad\quad\quad\quad\quad\quad
\quad\quad\quad\quad\quad \quad\quad(2.3)
\end{align*}

On the other hand, we see that
\begin{align*}
\|f\chi_{j}\|_{L^{p(\cdot)}(\omega)}&=2^{-j\alpha}\left(2^{j\alpha q}\|f\chi_{j}\|_{L^{p(\cdot)}(\omega)}^{q}\right)^{1/q}\\
&\lesssim 2^{-j\alpha}\left(\sum_{l=-\infty}^{j}2^{l\alpha q}\|f\chi_{l}\|_{L^{p(\cdot)}(\omega)}^{q}\right)^{1/q}\\
&=C2^{j(\lambda-\alpha)}\left(2^{-j\lambda}\left(\sum_{l=-\infty}^{j}2^{l \alpha q}\|f\chi_{l}\|_{L^{p(\cdot)}(\omega)}^{q}\right)^{1/q}\right)\\
&\lesssim 2^{j(\lambda-\alpha)}\|f\|_{M\dot{K}^{\alpha,\lambda}_{q,p(\cdot)}(\omega)}.\quad\quad\quad\quad\quad\quad\quad\quad
\quad\quad\quad\quad\quad \quad\quad\quad\quad\quad\quad\quad(2.4)
\end{align*}

This implies that in view of $n(1-r)+\lambda-\alpha>0$, combing  (2.3) and (2.4), we get

\begin{align*}
I_{1}&\lesssim m^{\frac{n}{2}} \cdot \sup_{L\in \mathbb{Z}}2^{-L\lambda q}\sum_{k=-\infty}^{L}2^{\alpha q k}\left(\sum_{j=-\infty}^{k-2}2^{(j-k)n(1-r)}\|f\chi_{j}\|_{L^{p(\cdot)}(\omega)}\right)^{q}\\
&\lesssim m^{\frac{n}{2}}  \|f\|_{M\dot{K}^{\alpha,\lambda}_{q,p(\cdot)}(\omega)}^{q}\sup_{L\in \mathbb{Z}}2^{-L\lambda q}\sum_{k=-\infty}^{L}2^{\lambda q k}\left(\sum_{j=-\infty}^{k-2}2^{(j-k)(n(1-r)+\lambda-\alpha)}\right)\\
&\lesssim m^{\frac{n}{2}}  \|f\|_{M\dot{K}^{\alpha,\lambda}_{q,p(\cdot)}(\omega)}^{q}\sup_{L\in \mathbb{Z}}2^{-L\lambda q}\left(\sum_{k=-\infty}^{L}2^{\lambda q k}\right)^{q}\\
&\lesssim m^{\frac{n}{2}}  \|f\|_{M\dot{K}^{\alpha,\lambda}_{q,p(\cdot)}(\omega)}^{q}.
\end{align*}
For $I_{2}$, according to Lemma 2.6, we get the following estimate
$$
I_{2}\lesssim m^{\frac{n}{2}} \sup_{L\in \mathbb{Z}}2^{-L\lambda q}\sum_{k=-\infty}^{L}2^{\alpha q k}\|f\chi_{k}\|_{L^{p(\cdot)}(\omega)}^{q}\lesssim m^{\frac{n}{2}} \|f\|_{M\dot{K}^{\alpha,\lambda}_{q,p(\cdot)}(\omega)}^{q}.
$$

We now turn to estimate $I_{3}$. Similar to $I_{1}$, we deduce that
$$
|T_{m,l}(f_{j})(x)|=|\int_{\mathbb{R}^{n}}K(x,y)f_{j}(y)\mathrm{d}y|\leq C\int_{\mathbb{R}^{n}}\frac{1}{|x-y|^{n}}f_{j}(y)\mathrm{d}y \approx C 2^{-jn}\|f\|_{L^{1}(\mathbb{R}^{n})}.
$$
 As argued before, by applying Lemmas 2.1 and 2.3, we obtain
\begin{align*}
|T_{m,l}(f_{j})\chi _{k}&\leq C 2^{-jn}\|f\omega^{\frac{1}{p(\cdot)}}\chi_{j}\|_{L^{p(\cdot)}(\mathbb{R}^{n})}\|f\omega^{-\frac{1}{p(\cdot)}}\chi_{j}\|_{L^{p^{\prime}(\cdot)}(\mathbb{R}^{n})}\\
&\leq C\|f\chi_{j}\|_{L^{p(\cdot)}(\omega)}\|\chi_{B_{j}}\|_{L^{p(\cdot)}(\omega)}^{-1}.\quad\quad\quad\quad\quad\quad\quad\quad
\quad\quad\quad\quad\quad \quad\quad(2.5)
\end{align*}

By using Lemma 2.2, it follows
\begin{align*}
\|T_{m,l}(f_{j})\chi_{k}\|_{L^{p(\cdot)}(\omega)}&\leq C\|f\chi_{j}\|_{L^{p(\cdot)}(\omega)}\|\chi_{B_{j}}\|_{L^{p(\cdot)}(\omega)}^{-1}\|\chi_{B_{k}}\|_{L^{p(\cdot)}(\omega)}\\
&\leq C 2^{(k-j)n\delta}\|f\chi_{j}\|_{L^{p(\cdot)}(\omega)}.\quad\quad\quad\quad\quad\quad\quad
\quad\quad\quad\quad\quad \quad\quad(2.6)
\end{align*}

From (2.4) and (2.6), in view of $\alpha+n\delta-\lambda>0$, we conclude that
\begin{align*}
I_{3}&\lesssim m^{\frac{n}{2}} \cdot \sup_{L\in \mathbb{Z}}2^{-L\lambda q}\sum_{k=-\infty}^{L}2^{\alpha q k}\left(\sum_{j=k+2}^{\infty}2^{(k-j)n\delta}\|f\chi_{j}\|_{L^{p(\cdot)}(\omega)}\right)^{q}\\
&\lesssim m^{\frac{n}{2}} \cdot \|f\|_{M\dot{K}^{\alpha,\lambda}_{q,p(\cdot)}(\omega)}^{q}\sup_{L\in \mathbb{Z}}2^{-L\lambda q}\sum_{k=-\infty}^{L}2^{\lambda q k}\left(\sum_{j=k+2}^{\infty}2^{(k-j)(n\delta+\alpha-\lambda)}\right)^{q}\\
&\lesssim m^{\frac{n}{2}} \cdot \|f\|_{M\dot{K}^{\alpha,\lambda}_{q,p(\cdot)}(\omega)}^{q}\sup_{L\in \mathbb{Z}}2^{-L\lambda q}\left(\sum_{k=-\infty}^{L}2^{\lambda q k}\right)\\
&\lesssim m^{\frac{n}{2}} \cdot \|f\|_{M\dot{K}^{\alpha,\lambda}_{q,p(\cdot)}(\omega)}^{q}.
\end{align*}

Combing the estimates of $I_{1}$, $I_{2}$ and $I_{3}$, we finish the proof of \textbf{Lemma 2.9}.

\section*{\normalsize\HT \ 3~~~Proofs of Theorems 1.1-1.5}

We begin with the following preliminary lemmas.\\
$\textbf{{Lemmas}~3.1}$ Let $0\leq\lambda<\infty,~p(\cdot)\in \mathcal{P}(\mathbb{R}^{n})\cap  LH(\mathbb{R}^{n}),~0<q<\infty,~\omega\in A_{rp(\cdot)}$, $1/p_{-}<r<1$  and $-n\delta+\lambda<\alpha<n(1-r)+\lambda$, where $0<\delta<1$ is the constant appearing in Lemma 2.2. and $t(x)$ be a homogeneous of degree $-n-1$ and locally integrable in $|x|>0$. Let $b\in Lip(\mathbb{R}^{n})$  and $K$ is defined by
$$
Kf(x)=\lim_{\varepsilon\mapsto0}\int_{|x-y|>\varepsilon}t(x)(b(x)-b(y))f(y)\mathrm{d}y.
$$
If $s(x)\in\mathcal{C}^{1}(S^{n-1})$ and $\int_{S^{n-1}}s(x)x_{j}\mathrm{d}\sigma(x)=0,~j=1,\ldots,n,$  then
$$~\|Kf\|_{M\dot{K}^{\alpha,\lambda}_{q,p(\cdot)}(\omega)}\lesssim(\|\nabla s\|_{L^{\infty}(S^{n-1})}+\|s\|_{L^{\infty}(S^{n-1})})\|\nabla b\|_{L^{\infty}}\|f\|_{M\dot{K}^{\alpha,\lambda}_{q,p(\cdot)}(\omega)}.$$

\textbf{Proof of Lemma 3.1}.

Let $k(x,y)=s(x-y)(b(x)-b(y))$, the kernel k satisfies the following inequalities for all $x,x_{0},y\in\mathbb{R}^{n}$ with $|x-x_{0}|\leq 1/2 |y-x|$,
$$
|k(x,y)-k(x_{0},y)|\lesssim\|\nabla s\|_{L^{\infty}(S^{n-1})}\|\nabla b\|_{L^{\infty}(S^{n-1})}|x-x_{0}||y-x|^{-n-1}. \eqno(3.1)
$$
and
$$
|k(x,y)|\lesssim\|s\|_{L^{\infty}(S^{n-1})}\|\nabla b\|_{L^{\infty}(S^{n-1})}|y-x|^{-n}. \eqno(3.2)
$$
This, together with the boundedness of $K$ on $L^{2}(\mathbb{R}^{n})$(see [41]), tells us $K$ is a generalized Calder\'{o}n-Zygmund operator and is bounded on $M\dot{K}^{\alpha,\lambda}_{q,p(\cdot)}(\omega)$ with bound $(\|\nabla s\|_{L^{\infty}}+\|s\|_{L^{\infty}})\|\nabla b\|_{L^{\infty}}\|f\|_{M\dot{K}^{\alpha,\lambda}_{q,p(\cdot)}(\omega)}$
(see [40]).

Thus, we get
$$
\|Kf\|_{M\dot{K}^{\alpha,\lambda}_{q,p(\cdot)}(\omega)}\lesssim(\|\nabla s\|_{L^{\infty}(S^{n-1})}+\|s\|_{L^{\infty}(S^{n-1})})\|\nabla b\|_{M\dot{K}^{\alpha,\lambda}_{q,p(\cdot)}(\omega)}.
$$
This completes our proof.

$\textbf{{Lemmas}~3.2}$ Let $0\leq\lambda<\infty,~p(\cdot)\in \mathcal{P}(\mathbb{R}^{n})\cap  LH(\mathbb{R}^{n}),~0<q<\infty,~\omega\in A_{rp(\cdot)}$, $1/p_{-}<r<1$  and $-n\delta+\lambda<\alpha<n(1-r)+\lambda$, where $0<\delta<1$ is the constant appearing in Lemma 2.2. Suppose that $b\in Lip(\mathbb{R}^{n})$ and $S$  be a singular operator which is defined by
$$
Sf(x)=\lim_{\varepsilon\mapsto0}\int_{|x-y|>\varepsilon}K(x-y)f(y)\mathrm{d}y,
$$
where $K(x)\in\mathcal{C}^{3}(S^{n-1})$, $\int_{S^{n-1}}K(x)\mathrm{d}\sigma(x)=0$ and $K(\lambda x)=\lambda^{-n}K(x)$, for $x\in\mathbb{R}^{n}\backslash\{0\}$,\\~$\lambda>0$, then we have, for $f\in C_{0}^{\infty}$,
$$
\|[b,S]\frac{\partial f}{\partial x_{j}}\|_{M\dot{K}^{\alpha,\lambda}_{q,p(\cdot)}(\omega)}\lesssim\max_{|\beta|\leq2}\|\partial^{\beta}K\|_{L^{\infty}(S^{n-1})}\|\nabla b\|_{L^{\infty}}\|f\|_{M\dot{K}^{\alpha,\lambda}_{q,p(\cdot)}(\omega)}.
$$

\textbf{Proof of Lemma 3.2}.
We remark that $f\in C_{0}^{\infty}$  leads to  $f \in L^{p(\cdot)}(\omega)([42])$. Let $f_{j},~b_{j},~K_{j}$ stands for $\frac{\partial f}{\partial x_{j}},~\frac{\partial b}{\partial x_{j}},~\frac{\partial K}{\partial x_{j}}$, respectively. Write
$$
[b,S]f_{j}(x)=\lim_{\varepsilon\mapsto0}\int_{|x-y|>\varepsilon}K(x-y)(b(x)-b(y))f_{j}(y)\mathrm{d}y.
$$
With the way similar to that used in [33, Lemma 5.2], it is not difficult to obtain Lemma 3.2.  Hence we omit details here.

\textbf{Proof of Theorem 1.1}.

Let
$$
\Omega(x,y)=\sum_{m\geq1}\sum^{d_{m}}_{j=1}a_{m,j}(x)Y_{m,j}(y).
$$
From [32],  for any $x$, we can write the coefficients $a_{m,j}$ as
$$
a_{m,j}(x)=(-1)^{n}m^{-n}(m+n-2)^{-n}\int_{S^{n-1}}L_{y^{'}}^{n}(\Omega(x,y^{'}))Y_{m,j}(y^{'})\mathrm{d}\lambda(y^{'}),m\geq1, \eqno(3.3)
$$
where $L(F)=|x|^{2}\Delta F(x).$

We will firstly prove the conclusion  {\it{(1)}}.  Write
\begin{align*}
(TD^{\gamma}-D^{\gamma}T)f & =\sum_{m=1}^{\infty}\sum_{j=1}^{d_{m}}(a_{m,j}T_{m,j}D^{\gamma}-D^{\gamma}a_{m,j}T_{m,j})f\\
&=\sum_{m=1}^{\infty}\sum_{j=1}^{d_{m}}(a_{m,j}D^{\gamma}T_{m,j}-D^{\gamma}a_{m,j}T_{m,j})f\\
&=\sum_{m=1}^{\infty}\sum_{j=1}^{d_{m}}[a_{m,j},D^{\gamma}]T_{m,j}f.
\end{align*}
By condition (3.3), it follows that
$$
D^{\gamma}a_{m,j}(x)=(-1)^{n}m^{-n}(m+n-2)^{-n}\int_{S^{n-1}}D_{x}^{\gamma}L_{y^{'}}^{n}(\Omega(x,y^{'}))Y_{m,j}(y^{'})\mathrm{d}\lambda(y^{'}), m\geq1.
$$
Further, by applying the condition (1.4), we have
$$
\|D^{\gamma}a_{m,j}\|_{L^{\infty}}\lesssim m^{-2n}. \eqno(3.4)
$$
Morever, By the fact that $[b,D^{\gamma}]$ is a generalized Calder\'{o}n-Zygmund operator (see [42]), which is defined by
$$
[b,D^{\gamma}]f(x)=C(\gamma)\int_{\mathbb{R}^{n}}\frac{(b(x)-b(y))}{|x-y|^{n+\gamma}}f(y)\mathrm{d}y.
$$
Thus, we see that $[b,D^{\gamma}]f(x)$ is bounded on the  $M\dot{K}^{\alpha,\lambda}_{q,p(\cdot)}(\omega)$ by applying Lemma 2.8. Namely
$$
\|[b,D^{\gamma}]f\|_{M\dot{K}^{\alpha,\lambda}_{q,p(\cdot)}(\omega)}\lesssim \|D^{^{\gamma}}b\|_{BMO}\|f\|_{M\dot{K}^{\alpha,\lambda}_{q,p(\cdot)}(\omega)}. \eqno(3.5)
$$
Then by $d_{m}\simeq m^{n-2}$ (see [30]),  (3.4), (3.5) and Lemma 2.9,  we have
 \begin{align*}
\|(TD^{\gamma}-D^{\gamma}T)f\|_{M\dot{K}^{\alpha(\cdot),q}_{p(\cdot),\lambda}}&\leq\sum_{m=1}^{\infty}\sum_{j=1}^{d_{m}}
\|[a_{m,j},D^{\gamma}]T_{m,j}f\|_{{M\dot{K}^{\alpha,\lambda}_{q,p(\cdot)}(\omega)}}\\
&\lesssim \sum_{m=1}^{\infty}\sum_{j=1}^{d_{m}}\|D^{\gamma}a_{m,j}\|_{BMO}\|T_{m,j}f\|_{_{M\dot{K}^{\alpha,\lambda}_{q,p(\cdot)}(\omega)}}\\
&\lesssim \sum_{m=1}^{\infty}\sum_{j=1}^{d_{m}}m^{\frac{n}{2}}\|D^{\gamma}a_{m,j}\|_{L^{\infty}}\|f\|_{_
{M\dot{K}^{\alpha,\lambda}_{q,p(\cdot)}(\omega)}}\\
&\lesssim \sum_{m=1}^{\infty}m^{n-2}m^{\frac{n}{2}}m^{-2n}\|f\|_{_{M\dot{K}^{\alpha,\lambda}_{q,p(\cdot)}(\omega)}}\\
&\lesssim\|f\|_{_{M\dot{K}^{\alpha,\lambda}_{q,p(\cdot)}(\omega)}}.
\end{align*}

Now let us turn to estimate {\it{(2)}}. By applying the definition of $T^{\sharp}$ and $T^{\ast}$ we can deduce that
$$
 (T^{\sharp}-T^{\ast})D^{\gamma}f=\sum_{m=1}^{\infty}\sum_{j=1}^{d_{m}}(-1)^{m}[\bar{a}_{m,j},T_{m,j}]D^{\gamma}f. \eqno(3.6)
$$
To estimate $M\dot{K}^{\alpha,\lambda}_{q,p(\cdot)}(\omega)$ norm of $(T^{\ast}-T^{\sharp})D^{\gamma}$, we first consider
$[b,T_{m.j}]D^{\gamma}$ for any fixed $b\in I_{\gamma}(\mathrm{BMO})$. Noting that $b(x)-b(y)=(b(x)-b(z))-(b(y)-b(z))$, for any $x,y,z\in \mathbb{R}^{n}$,  then we have
$$
[b,T_{m,j}]D^{\gamma}f=[b,D^{\gamma}T_{m,j}]f-T_{m,j}[b,D^{\gamma}]f.
$$
Thus,  we get by (3.5) and Lemma 2.9
$$
 \|T_{m,j}[b,D^{\gamma}]f\|_{{M\dot{K}^{\alpha,\lambda}_{q,p(\cdot)}(\omega)}}\lesssim m^{\frac{n}{2}}\|D^{\gamma}b\|_{BMO}\|f\|_{_{M\dot{K}^{\alpha,\lambda}_{q,p(\cdot)}(\omega)}} .\eqno(3.7)
$$
Further, we estimate the $M\dot{K}^{\alpha,\lambda}_{q,p(\cdot)}(\omega)$ norm of $[b,D^{\gamma}T_{m,j}]f$.  From the fact that $[b,D^{\gamma}T_{m,j}]$ is a generalized Calder\'{o}n-Zygmund operator with kernel (see [33])
$$
|k_{m,j}(x,y)|\lesssim m^{\frac{n}{2}-1+\gamma}\|D^{\gamma}b\|_{BMO}\frac{1}{|x-y|^{n}},
$$
then we get by Lemma 2.8
$$
  \|[b,D^{\gamma}T_{m,j}]f\|_{_{M\dot{K}^{\alpha,\lambda}_{q,p(\cdot)}(\omega)}}\lesssim m^{\frac{n}{2}+\gamma}\|D^{\gamma}b\|_{BMO}\|f\|_{_{M\dot{K}^{\alpha,\lambda}_{q,p(\cdot)}(\omega)}}.\eqno(3.8)
$$
Then, combining (3.7) with (3.8), we have
\begin{eqnarray*}
\|[b,T_{m,j}]D^{\gamma}f\|_{_{M\dot{K}^{\alpha,\lambda}_{q,p(\cdot)}(\omega)}} &&\lesssim m^{\frac{n}{2}+\gamma}\|D^{\gamma}b\|_{BMO}\|f\|_{_{M\dot{K}^{\alpha,\lambda}_{q,p(\cdot)}(\omega)}}\\
&&\quad+m^{\frac{n}{2}
}\|D^{\gamma}b\|_{BMO}\|f\|_{_{M\dot{K}^{\alpha,\lambda}_{q,p(\cdot)}(\omega)}}\\
&&\lesssim m^{\frac{n}{2}+\gamma}\|D^{\gamma}b\|_{BMO}\|f\|_{_{M\dot{K}^{\alpha,\lambda}_{q,p(\cdot)}(\omega)}}.~~ \quad\quad\quad\quad\quad\quad\quad\quad\quad\quad
\quad(3.9)
\end{eqnarray*}
By condition (3.4), (3.6) and (3.9), we get
\begin{align*}
\|(T^{\sharp}-T^{\ast})D^{\gamma}f\|_{_{M\dot{K}^{\alpha,\lambda}_{q,p(\cdot)}(\omega)}}&\leq \sum_{m=1}^{\infty}\sum_{j=1}^{d_{m}}\|[\bar{a}_{m,j},T_{m,j}]D^{\gamma}f\|_{_{M\dot{K}^{\alpha,\lambda}_{q,p(\cdot)}(\omega)}}\\
&\lesssim \sum_{m=1}^{\infty}\sum_{j=1}^{d_{m}}m^{\frac{n}{2}+\gamma}\|D^{\gamma}\bar{a}_{m,j}\|_{BMO}\|f\|_{_{M\dot{K}^{\alpha,\lambda}_{q,p(\cdot)}(\omega)}}\\
 &\lesssim\sum_{m=1}^{\infty}\sum_{j=1}^{d_{m}}m^{\frac{n}{2}+\gamma}\|D^{\gamma}\bar{a}_{m,j}\|_{L^{\infty}}\|f\|_{_{M\dot{K}^{\alpha,\lambda}_{q,p(\cdot)}(\omega)}}\\
&\lesssim\sum_{m=1}^{\infty}m^{n-2}m^{\frac{n}{2}+\gamma}m^{-2n}\|f\|_{_{M\dot{K}^{\alpha,\lambda}_{q,p(\cdot)}(\omega)}}\\
&\lesssim\|f\|_{_{M\dot{K}^{\alpha,\lambda}_{q,p(\cdot)}(\omega)}}.
\end{align*}
Thus we finish the proof of \textbf{Theorem 1.1}.

\textbf{Proof of Theorem 1.2}.
Let
$$
T_{1}f(x)= \int_{\mathbb{R}^{n}}\frac{\Omega_{1}(x,x-y)}{|x-y|^{n}}f(y)\mathrm{d}y ~~{\rm and}~~T_{2}f(x)= \int_{\mathbb{R}^{n}}\frac{\Omega_{2}(x,x-y)}{|x-y|^{n}}f(y)\mathrm{d}y.
$$
Write
$$
\Omega_{1}(x,y)=\sum_{m\geq1}\sum^{d_{m}}_{j=1}a_{m,j}(x)Y_{m,j}(y) ~~{\rm and}~~\Omega_{2}(x,y)=\sum_{\lambda\geq1}\sum^{d_{\lambda}}_{\mu=1}b_{\lambda,\mu}(x)Y_{\lambda,\mu}(y),
$$
where
 $$
a_{m,j}(x)=\int_{S^{n-1}}\Omega_1(x,z^{\prime})\overline{Y_{m,j}(z^{\prime})}\mathrm{d}\lambda(z^{\prime})
 ~~{\rm and}~~b_{\lambda,\mu}(x)=\int_{S^{n-1}}\Omega_{2}(x,z^{\prime})\overline{Y_{\lambda,\mu}(z^{\prime})}\mathrm{d}\lambda(z^{\prime}).
$$

For any $x\in\mathbb{R}^{n}$, with a similar argument used in the proof of Theorem 1.1 in terms of (1.4) and (1.5), we can obtain that
$$
\|a_{m,j}\|_{L^{\infty}}\lesssim m^{-2n}. \eqno(3.10)
$$
$$
\|D^{\gamma}b_{\lambda,\mu}\|_{L^{\infty}}\lesssim m^{-2n}. \eqno(3.11)
$$

Let
$$
T_{m,j}f(x)=\frac{Y_{m,j}}{|\cdot|^{n}}\ast f(x) ~~{\rm and~~}
T_{\lambda,\mu}f(x)=\frac{Y_{\lambda,\mu}}{|\cdot|^{n}}\ast f(x).
$$
Since $\Omega_{1}(x,y)$ and $\Omega_{2}(x,y)$ satisfy (1.3), then we get
 $$T_{1}f(x)=\sum_{m\geq 1}\sum_{j=1}^{d_{m}}a_{m,j}(x)T_{m,j}f(x)~~{\rm and}~~T_{2}f(x)=\sum_{\lambda\geq 1}\sum_{\mu=1}^{d_{\lambda}}b_{\lambda,\mu}(x)T_{\lambda,\mu}f(x).$$
 Write(see [33])
 $$
 (T_{1}\circ T_{2})f(x)=\sum_{m=1}^{\infty}\sum_{j=1}^{d_{m}}\sum_{\lambda=1}^{\infty}\sum_{\mu=1}^{d_{\lambda}}a_{m,j}(x)b_
 {\lambda,\mu}(x)(T_{m,j}T_{\lambda,\mu}f)(x),
 $$
 $$
 (T_{1}T_{2})f(x)=\sum_{m=1}^{\infty}\sum_{j=1}^{d_{m}}\sum_{\lambda=1}^{\infty}\sum_{\mu=1}^{d_{\lambda}}a_{m,j}T_{m,j}
 (b_{\lambda,\mu}T_{\lambda,\mu}f)(x).
 $$
 Then
\begin{align*}
  (T_{1}\circ T_{2}-T_{1}T_{2})D^{\gamma}&f=\sum_{m=1}^{\infty}\sum_{j=1}^{d_{m}}\sum_{\lambda=1}^{\infty}\sum_{\mu=1}
 ^{d_{\lambda}}a_{m,j}(b_{\lambda,\mu}(x)T_{m,j}-T_{m,j}b_{\lambda,\mu}(x))T_{\lambda,\mu}D^{\gamma}f\\
  &=\sum_{m=1}^{\infty}\sum_{j=1}^{d_{m}}\sum_{\lambda=1}^{\infty}\sum_{\mu=1}
 ^{d_{\lambda}}a_{m,j}(b_{\lambda,\mu}(x)T_{m,j}-T_{m,j}b_{\lambda,\mu}(x))D^{\gamma}T_{\lambda,\mu}f\\
 &=\sum_{m=1}^{\infty}\sum_{j=1}^{d_{m}}\sum_{\lambda=1}^{\infty}\sum_{\mu=1}
 ^{d_{\lambda}}a_{m,j}[b_{\lambda,\mu},T_{m,j}]D^{\gamma}T_{\lambda,\mu}f.
\end{align*}
Therefore, together with (3.9), (3.10), (3.11) and Lemma 2.9, we obtain
\begin{align*}
&\|(T_{1}\circ T_{2}-T_{1}T_{2})D^{\gamma}f\|_{_{M\dot{K}^{\alpha,\lambda}_{q,p(\cdot)}(\omega)}}\\
&\lesssim\sum_{m=1}^{\infty}\sum_{j=1}^{d_{m}}\sum_{\lambda=1}^{\infty}\sum_{\mu=1}^{d_{\lambda}}\|a_{m,j}\|_
{L^{\infty}}\|[b_{\lambda,\mu},T_{m,j}]D^{\gamma}T_{\lambda,\mu}f\|_{_{M\dot{K}^{\alpha,\lambda}_{q,p(\cdot)}(\omega)}}\\
&\lesssim\sum_{m=1}^{\infty}\sum_{j=1}^{d_{m}}\sum_{\lambda=1}^{\infty}\sum_{\mu=1}^{d_{\lambda}}\|a_{m,j}\|_
{L^{\infty}}\|D^{\gamma}b_{\lambda,\mu}\|_{BMO} m^{\frac{n}{2}+\gamma}\|T_{\lambda,\mu}f\|_{_{M\dot{K}^{\alpha,\lambda}_{q,p(\cdot)}(\omega)}}\\
&\lesssim\sum_{m=1}^{\infty}\sum_{j=1}^{d_{m}}\sum_{\lambda=1}^{\infty}\sum_{\mu=1}^{d_{\lambda}}\|a_{m,j}\|_
{L^{\infty}}\|D^{\gamma}b_{\lambda,\mu}\|_{L^{\infty}}m^{\frac{n}{2}+\gamma}\lambda^{\frac{n}{2}}\|f\|_{_{M\dot{K}^{\alpha,\lambda}_{q,p(\cdot)}(\omega)}}\\
&\lesssim \sum_{m=1}^{\infty}m^{n-2}m^{-2n}m^{\frac{n}{2}+\gamma}\sum_{\lambda=1}^{\infty}\lambda^{n-2}\lambda^{-2n}\lambda
^{\frac{n}{2}}\|f\|_{_{M\dot{K}^{\alpha,\lambda}_{q,p(\cdot)}(\omega)}}\\
&\lesssim\|f\|_{_{M\dot{K}^{\alpha,\lambda}_{q,p(\cdot)}(\omega)}}.
\end{align*}
This finishes the proof of \textbf{Theorem 1.2}.

\textbf{Proof of Theorem 1.3}.
We estimate that term exactly as we did for the corresponding boundedness in Theorem 1.1 in the above arguments. Without loss of generality, we only to prove  {\it{(2)}} and  {\it{(3)}} of Theorem 1.3. By using the fact that $\Omega_{1}(x,y)$ and $\Omega_{2}(x,y)$ satisfy (1.5), therefore, we have shown that
  $$
\|a_{m,j}\|_{L^{\infty}}\lesssim m^{-2n},~~~\|b_{\lambda,\mu}\|_{L^{\infty}}\lesssim \lambda^{-2n}. \eqno(3.12)
$$
Firstly, let's prove  {\it{(2)}}. As in the proof of Theorem 1.1, we can get
$$
(T_{1}^{\sharp}-T_{1}^{\ast})\mathcal{I}f=\sum_{m=1}^{\infty}\sum_{j=1}^{d_{m}}(-1)^{m}[\bar{a}_{m,j},T_{m,j}]\mathcal{I}f.
$$
we showed that $[b,T_{m,j}]$ is a special Calder\'{o}n-Zygmund operator, so it is a bounded operator on the $M\dot{K}^{\alpha,\lambda}_{q,p(\cdot)}(\omega)$ by applying Lemma 2.8. Thus we have
  $$
\|[b,T_{m,j}]f\|_{M\dot{K}^{\alpha,\lambda}_{q,p(\cdot)}(\omega)}\lesssim m^{\frac{n}{2}}\|b\|_{L^{\infty}}\|f\|_{M\dot{K}^{\alpha,\lambda}_{q,p(\cdot)}(\omega)}. \eqno(3.13)
$$
Then by (3.12), we get
\begin{align*}
\|(T_{1}^{\sharp}-T_{1}^{\ast})\mathcal{I}f\|_{M\dot{K}^{\alpha,\lambda}_{q,p(\cdot)}(\omega)}&\lesssim\sum_{m=1}^{\infty}m^{n-2}
m^{-3n/2}\|f\|_{M\dot{K}^{\alpha,\lambda}_{q,p(\cdot)}(\omega)}\\&\lesssim\|f\|_{M\dot{K}^{\alpha,\lambda}_{q,p(\cdot)}(\omega)}.
\end{align*}
Thus the conclusion {\it{(2)}} is proved. We now estimate {\it{(3)}}.  Write
$$
(T_{1}\circ T_{2}-T_{1}T_{2})\mathcal{I}f=\sum_{m=1}^{\infty}\sum_{j=1}^{d_{m}}\sum_{\lambda=1}^{\infty}\sum_{\mu=1}^{d_{\lambda}}[b_{\lambda,\mu},T_{m,j}]T_{\lambda,\mu}\mathcal{I}f.
$$
Therefore, by  (3.12), (3.13) and Lemma 2.9, we get
\begin{align*}
&\|(T_{1}\circ T_{2}-T_{1}T_{2})\mathcal{I}f\|_{M\dot{K}^{\alpha,\lambda}_{q,p(\cdot)}(\omega)}\\
&\lesssim\sum_{m=1}^{\infty}\sum_{j=1}^{d_{m}}\sum_{\lambda=1}^{\infty}\sum_{\mu=1}^{d_{\lambda}}\|a_{m,j}\|_{L^{\infty}}\|b_{\lambda,\mu}\|_{L^{\infty}}m^{\frac{n}{2}}\|T_{\lambda,\mu}\mathcal{I}f]\|
_{M\dot{K}^{\alpha,\lambda}_{q,p(\cdot)}(\omega)}.\\
&\lesssim\sum_{m=1}^{\infty}m^{n-2}m^{-2n}m^{n/2}\sum_{\lambda=1}^{\infty}\lambda^{n-2}\lambda^{-2n}\lambda^{n/2}\|f\|_{M\dot{K}^{\alpha,\lambda}_{q,p(\cdot)}(\omega)}\\
&\lesssim \|f\|_{M\dot{K}^{\alpha,\lambda}_{q,p(\cdot)}(\omega)}.
\end{align*}
Thus the conclusion {\it{(3)}} is also proved. Hence the proof of \textbf{Theorem 1.3} is finished.

\textbf{Proof of Theorem 1.4}.
In the first place, we will prove the conclusion (1). Write $D=\sum_{k=1}^{n}\mathcal{R}_{k}\frac{\partial}{\partial x_{k}}$, where $\mathcal{R}_{k}$ denotes the Riesz transform. As in the proof of Theorem 1.1, we have
\begin{align*}
&(TD-DT)f(x)=\sum_{m=1}^{\infty}\sum_{j=1}^{d_{m}}[a_{m,j},D]T_{m,j}f(x)\\
&=\sum_{m=1}^{\infty}\sum_{j=1}^{d_{m}}\sum_{k=1}^{n}\mathcal{R}_{k}[a_{m,j},\frac{\partial}{\partial x_{k}}]T_{m,j}f(x)+\sum_{m=1}^{\infty}\sum_{j=1}^{d_{m}}\sum_{k=1}^{n}[a_{m,j},\mathcal{R}_{k}]\frac{\partial}{\partial x_{k}}(T_{m,j}f)(x)\\
&=: J_{1}+J_{2}.
\end{align*}
We have by the Leibniz's rules that
$$
J_{1}=\sum_{m=1}^{\infty}\sum_{j=1}^{d_{m}}\sum_{k=1}^{n}\mathcal{R}_{k}(\frac{\partial}{\partial x_{k}}(a_{m,j})T_{m,j}f).
$$
Thus we deduce from (3.3) that
$$
\frac{a_{m,j}}{\partial x_{k}}(x)=(-1)^{n}m^{-n}(m+n-2)^{-n}\int_{S^{n-1}}\partial_{x_{k}}L_{y^{'}}^{n}(\Omega(x,y^{'}))Y_{m,j}(y^{'})\mathrm{d}\lambda(y^{'}),m\geq1.
$$
From this and (1.6), we get for $k=1,\ldots,n,$
 $$
\left\|\frac{\partial a_{m,j}}{\partial x_{k}}\right\|_{L^{\infty}}\lesssim m^{-2n}. \eqno(3.14)
$$
By using the fact that $\|\mathcal{R}_{k}g\|_{WM\dot{K}^{\alpha,\lambda}_{p,1}(\mathbb{R}^{n})}\lesssim\|g\|_{M\dot{K}^{\alpha,\lambda}_{p,1}(\mathbb{R}^{n})}$, $d_{m}\simeq m^{n-2}$ and Lemma 3.1,  then we have
\begin{align*}
\|J_{1}\|_{M\dot{K}^{\alpha,\lambda}_{q,p(\cdot)}(\omega)}&\lesssim\sum_{m=1}^{\infty}\sum_{j=1}^{d_{m}}\sum_{k=1}^{n}\|\mathcal{R}_{k}(\frac{\partial}{\partial x_{k}}(a_{m,j})T_{m,j}f)\|_{M\dot{K}^{\alpha,\lambda}_{q,p(\cdot)}(\omega)}\\
&\lesssim\sum_{m=1}^{\infty}\sum_{j=1}^{d_{m}}m^{-2n}m^{n/2}\|f\|_{M\dot{K}^{\alpha,\lambda}_{q,p(\cdot)}(\omega)}\\
&\lesssim\sum_{m=1}^{\infty}m^{n-2}m^{-2n}m^{n/2}\|f\|_{M\dot{K}^{\alpha,\lambda}_{q,p(\cdot)}(\omega)}\\
&\lesssim\|f\|_{M\dot{K}^{\alpha,\lambda}_{q,p(\cdot)}(\omega)}.
\end{align*}
By Lemma 3.2 and (3.14), a trivial computation shows that for $I_{2}$,
\begin{align*}
\|J_{2}\|_{M\dot{K}^{\alpha,\lambda}_{q,p(\cdot)}(\omega)}&\lesssim\sum_{m=1}^{\infty}\sum_{j=1}^{d_{m}}\sum_{k=1}^{n}\|\nabla a_{m,j}\|_{L^{\infty}}\|T_{m,j}f\|_{M\dot{K}^{\alpha,\lambda}_{q,p(\cdot)}(\omega)}\\
&\lesssim\sum_{m=1}^{\infty}\sum_{j=1}^{d_{m}}m^{-2n}m^{n/2}\|f\|_{M\dot{K}^{\alpha,\lambda}_{q,p(\cdot)}(\omega)}\\
&\lesssim\sum_{m=1}^{\infty}m^{n-2}m^{-2n}m^{n/2}\|f\|_{M\dot{K}^{\alpha,\lambda}_{q,p(\cdot)}(\omega)}\\
&\lesssim\|f\|_{M\dot{K}^{\alpha,\lambda}_{q,p(\cdot)}(\omega)}.
\end{align*}
Combining the estimates above, we arrive at the desired boundedness
$$
\|(TD-DT)f\|_{M\dot{K}^{\alpha,\lambda}_{q,p(\cdot)}(\omega)}\lesssim \|f\|_{M\dot{K}^{\alpha,\lambda}_{q,p(\cdot)}(\omega)}.
$$
We posterior prove the conclusion (2). Write $D=\sum_{k=1}^{n}\mathcal{R}_{k}\frac{\partial}{\partial x_{k}}$, we have
\begin{align*}
(T^{\sharp}-T^{\ast})Df(x)&=\sum_{m=1}^{\infty}\sum_{j=1}^{d_{m}}(-1)^{m}[\bar{a}_{m,j},T_{m,j}]Df(x)\\
&=\sum_{k=1}^{n}\sum_{m=1}^{\infty}\sum_{j=1}^{d_{m}}(-1)^{m}[\bar{a}_{m,j},T_{m,j}]\frac{\partial}{\partial x_{k}}(\mathcal{R}_{k}f)(x).\quad\quad\quad\quad\quad\quad\quad\quad\quad(3.15)
\end{align*}
We now turn to estimate the $M\dot{K}^{\alpha,\lambda}_{q,p(\cdot)}(\omega)$ norm of $ [\bar{a}_{m,j},T_{m,j}]\frac{\partial}{\partial x_{k}}(\mathcal{R}_{k}f)$. Applying (3.14), Lemma 3.2 and the fact that for any multi-index $\beta$
and $x\in \mathbb{R}^{n}\backslash\{0\},~m=1,2,\ldots.$ (see [32]),
$$
|\partial^{\beta}(|x|^{m})Y_{m,j}|\leq C(n)|x|^{m-|\beta|}m^{|\beta|+(n-2)/2}.\eqno(3.16)
$$
Hence, we get
\begin{align*}
\left\|[\bar{a}_{m,j},T_{m,j}]\frac{\partial}{\partial x_{k}}(\mathcal{R}_{k}f)\right\|_{M\dot{K}^{\alpha,\lambda}_{q,p(\cdot)}(\omega)}&
\lesssim\|\nabla\bar{a}_{m,j}\|_{L^{\infty}}\max_{|\beta|\leq2}\|\partial^{\beta}Y_{m,j}\|_{L^{\infty}(S^{n-1})}\|\mathcal{R}_{k}f\|_{M\dot{K}^{\alpha,\lambda}_{q,p(\cdot)}(\omega)}\\
&\lesssim m^{-2n}m^{n/2+1}\|f\|_{M\dot{K}^{\alpha,\lambda}_{q,p(\cdot)}(\omega)}\\
&\lesssim m^{-3n/2+1}\|f\|_{M\dot{K}^{\alpha,\lambda}_{q,p(\cdot)}(\omega)}.\quad\quad\quad\quad\quad
\quad\quad\quad\quad\quad(3.17)
\end{align*}
Combining the estimates of (4.15) with (4.17), we have
\begin{align*}
\|(T^{\sharp}-T^{\ast})Df\|_{M\dot{K}^{\alpha,\lambda}_{q,p(\cdot)}(\omega)}&\lesssim\sum_{m=1}^{\infty}m^{n-2}m^{-3n/2+1}\|f\|_{M\dot{K}^{\alpha,\lambda}_{q,p(\cdot)}(\omega)}\\
&\lesssim\|f\|_{M\dot{K}^{\alpha,\lambda}_{q,p(\cdot)}(\omega)}.
\end{align*}
Consequently, the proof of \textbf{Theorem 1.4} is completed.

\textbf{Proof of Theorem 1.5}.
Similar to the proof of Theorem 1.2,  we easily see that
$$
(T_{1}\circ T_{2}-T_{1}T_{2})Df=\sum_{m=1}^{\infty}\sum_{d=1}^{d_{m}}\sum_{\lambda=1}^{\infty}\sum_{\mu=1}^{d_{\lambda}}a_{m,j}[b_{\lambda,\mu},T_{m,j}]DT_{\lambda,\mu}f,
$$
where $a_{m,j}$ and $b_{\lambda,\mu}$ are same to occur in the proof of Theorem 1.2.  By (1.5) and (1.6), we have
 $$
\|a_{m,j}\|_{L^{\infty}}\lesssim m^{-2n}. \eqno(3.18)
$$
 $$
\|\nabla b_{\lambda,\mu}\|_{L^{\infty}}\lesssim \lambda^{-2n}. \eqno(3.19)
$$
write $D=\sum_{k=1}^{n}\frac{\partial}{\partial x_{k}}\mathcal{R}_{k},$ it then follows that
\begin{align*}
&\|(T_{1}\circ T_{2}-T_{1}T_{2})Df\|_{M\dot{K}^{\alpha,\lambda}_{q,p(\cdot)}(\omega)}\\
&\lesssim\sum_{m=1}^{\infty}\sum_{j=1}^{d_{m}}\sum_{\lambda=1}^{\infty}\sum_{\mu=1}^{d_{\lambda}}\|a_{m,j}\|_{L^{\infty}}\left\|[b_{\lambda,\mu},T_{m,j}](\sum_{k=1}^{n}\frac{\partial}{\partial x_{k}}\mathcal{R}_{k}T_{\lambda,\mu}f)\right\|_{M\dot{K}^{\alpha,\lambda}_{q,p(\cdot)}(\omega)}\\
&\lesssim\sum_{k=1}^{n}\sum_{m=1}^{\infty}\sum_{j=1}^{d_{m}}\sum_{\lambda=1}^{\infty}\sum_{\mu=1}^{d_{\lambda}}\|a_{m,j}\|_{L^{\infty}}\left\|[b_{\lambda,\mu},T_{m,j}](\frac{\partial}{\partial x_{k}}\mathcal{R}_{k}T_{\lambda,\mu}f)\right\|_{M\dot{K}^{\alpha,\lambda}_{q,p(\cdot)}(\omega)}.
\end{align*}
The above estimate, via Lemma 2.9, lead to that
\begin{align*}
&\|(T_{1}\circ T_{2}-T_{1}T_{2})Df\|_{M\dot{K}^{\alpha,\lambda}_{q,p(\cdot)}(\omega)}\\
&\lesssim\sum_{m=1}^{\infty}\sum_{j=1}^{d_{m}}\sum_{\lambda=1}^{\infty}\sum_{\mu=1}^{d_{\lambda}}\|a_{m,j}\|_{L^{\infty}}\|\nabla b_{\lambda,\mu}\|_{L^{\infty}}\max_{|\beta|\leq2}\|\partial^{\beta}Y_{m,j}\|_{L^{\infty}(S^{n-1})}\|T_{\lambda,\mu}\mathcal{R}_{k}f\|_{M\dot{K}^{\alpha,\lambda}_{q,p(\cdot)}(\omega)}.
\end{align*}
We thus obtain from (3.16), (3.18), (3.19) and Lemma 2.9 that,
\begin{align*}
&\|(T_{1}\circ T_{2}-T_{1}T_{2})Df\|_{M\dot{K}^{\alpha,\lambda}_{q,p(\cdot)}(\omega)}\\
&\lesssim\sum_{m=1}^{\infty}m^{n/2-1}m^{-2n}m^{n/2+1}\sum_{\lambda=1}^{\infty}\lambda^{n/2-1}\lambda^{-2n}\lambda^{n/2}\|f\|_{M\dot{K}^{\alpha,\lambda}_{q,p(\cdot)}(\omega)}\\
&\lesssim\|f\|_{M\dot{K}^{\alpha,\lambda}_{q,p(\cdot)}(\omega)}.
\end{align*}
Consequently, the proof of \textbf{Theorem 1.5} is finished.

By Remark 2.1, we can obtain the following outcomes.\\

$\textbf{{Corollary}~3.1}$  Let $p(\cdot)\in \mathcal{P}(\mathbb{R}^{n})\cap  LH(\mathbb{R}^{n}),~0<q<\infty,~\omega\in A_{rp(\cdot)}$, $1/p_{-}<r<1$  and $-n\delta<\alpha<n(1-r)$, where $0<\delta<1$ is the constant appearing in Lemma 2.2. Let $\Omega(x,y)$ which satisfies(1.2), (1.3) and (1.4), then we have

(1) $\|(TD^{\gamma}-D^{\gamma}T)f\|_{\dot{K}^{\alpha,q}_{p(\cdot)}(\omega)}\lesssim  \|f\|_{\dot{K}^{\alpha}_{q,p(\cdot)}(\omega)};$

(2) $\|(T^{\ast}-T^{\sharp})D^{\gamma}f\|_{\dot{K}^{\alpha,q}_{p(\cdot)}(\omega)}\lesssim \|f\|_{\dot{K}^{\alpha,q}_{p(\cdot)}(\omega)}.$\\
$\textbf{{Corollary}~3.2}$  Let $p(\cdot)\in \mathcal{P}(\mathbb{R}^{n})\cap  LH(\mathbb{R}^{n}),~0<q<\infty,~\omega\in A_{rp(\cdot)}$, $1/p_{-}<r<1$  and $-n\delta<\alpha<n(1-r)$, where $0<\delta<1$ is the constant appearing in Lemma 2.2. Suppose that $\Omega_{1}(x,y)$ and $\Omega_{2}(x,y)$  satisfy (1.2) and (1.3). If $\Omega_{2}(x,y)$  satisfies (1.4) and (1.5), then we have
$$
\|(T_{1}\circ T_{2}-T_{1}T_{2})D^{\gamma}f\|_{\dot{K}^{\alpha,q}_{p(\cdot)}(\omega)}\lesssim\|f\|_{\dot{K}^{\alpha,q}_{p(\cdot)}(\omega)}.
$$
$\textbf{{Corollary}~3.3}$  Let $p(\cdot)\in \mathcal{P}(\mathbb{R}^{n})\cap  LH(\mathbb{R}^{n}),~0<q<\infty,~\omega\in A_{rp(\cdot)}$, $1/p_{-}<r<1$  and $-n\delta<\alpha<n(1-r)$, where $0<\delta<1$ is the constant appearing in Lemma 2.2.  suppose that $\Omega_i(x,y)(i=1,2)$ satisfies(1.2),(1.3) and (1.5). Then we have

(1) $\|(T_{1}\mathcal{I}-\mathcal{I}T_{1})f\|_{\dot{K}^{\alpha,q}_{p(\cdot)}(\omega)}\lesssim  \|f\|_{\dot{K}^{\alpha,q}_{p(\cdot)}(\omega)};$

(2) $\|(T_{1}^{\ast}-T_{2}^{\sharp})\mathcal{I}f\|_{\dot{K}^{\alpha,q}_{p(\cdot)}(\omega)}\lesssim \|f\|_{\dot{K}^{\alpha,q}_{p(\cdot)}(\omega)};$

(3) $\|(T_{1}\circ T_{2}-T_{1}T_{2})\mathcal{I}f\|_{\dot{K}^{\alpha,q}_{p(\cdot)}(\omega)}\lesssim\|f\|_{\dot{K}^{\alpha,q}_{p(\cdot)}(\omega)}.$\\
$\textbf{{Corollary}~3.4}$  Let $p(\cdot)\in \mathcal{P}(\mathbb{R}^{n})\cap  LH(\mathbb{R}^{n}),~0<q<\infty,~\omega\in A_{rp(\cdot)}$, $1/p_{-}<r<1$  and $-n\delta<\alpha<n(1-r)$, where $0<\delta<1$ is the constant appearing in Lemma 2.2.  suppose that $\Omega(x,y)$ satisfies(1.2), (1.3) and (1.6), then we have

(1) $\|(TD-DT)f\|_{\dot{K}^{\alpha,q}_{p(\cdot)}(\omega)}\lesssim  \|f\|_{\dot{K}^{\alpha,q}_{p(\cdot)}(\omega)};$

(2) $\|(T^{\ast}-T^{\sharp})Df\|_{\dot{K}^{\alpha,q}_{p(\cdot)}(\omega)}\lesssim \|f\|_{\dot{K}^{\alpha,q}_{p(\cdot)}(\omega)}.$\\
$\textbf{{Corollary}~3.5}$  Let $p(\cdot)\in \mathcal{P}(\mathbb{R}^{n})\cap  LH(\mathbb{R}^{n}),~0<q<\infty,~\omega\in A_{rp(\cdot)}$, $1/p_{-}<r<1$  and $-n\delta<\alpha<n(1-r)$, where $0<\delta<1$ is the constant appearing in Lemma 2.2.  Suppose that $\Omega_{1}(x,y)$ and $\Omega_{2}(x,y)$ satisfies(1.2), (1.3). If $\Omega_{1}(x,y)$ satisfies (1.5) and $\Omega_{2}(x,y)$ satisfies (1.6). Then we have
$$
\|(T_{1}\circ T_{2}-T_{1}T_{2})Df\|_{\dot{K}^{\alpha,q}_{p(\cdot)}(\omega)}\lesssim\|f\|_{\dot{K}^{\alpha,q}_{p(\cdot)}(\omega)}.
$$

\noindent$\textbf{Acknowledgments}$ \\
This work is supported by the Doctoral Scientific Research Foundation of Northwest Normal University (202003101203) and Teachers Scientific Research Ability Promotion Project of Northwest Normal University (NWNU-LKQN2021-03), and Open Foundation of Hubei Key Laboratory of Applied Mathematics (Hubei University)(HBAM202205).


\def\refname{References }

\end{document}